
\documentstyle{amsppt}

\NoBlackBoxes
\magnification=\magstep1
\hsize=6.5truein
\vsize=8.5truein
\document

\topmatter 
\title 
The vector
$k$--constrained KP 
hierarchy and Sato's Grassmannian
\endtitle
\author
Johan van de Leur
\endauthor
\affil
Faculty of Applied Mathematics\\
University of Twente\\
P.O. Box 217\\
7500 AE Enschede\\
 The Netherlands\\
\endaffil

\thanks
The author is
supported by the ``Stichting Fundamenteel Onderzoek der
Materie (F.O.M.)''. E-mail: vdleur\@math.utwente.nl
\endthanks
\abstract
We use the representation theory of the infinite matrix group
to show that (in the polynomial case) the $n$--vector $k$--constrained
KP hierarchy has a natural geometrical interpretation on
Sato's infinite Grassmannian. This description generalizes the
the $k$--reduced KP or Gelfand--Dickey hierarchies. 

\endabstract
\endtopmatter

\break
\subheading{\S 1. Introduction}
\vskip 10pt
\noindent 
It is well--known that the $k$--th Gelfand--Dickey hierarchy, 
which generalizes the Korteweg--de Vries (KdV) hierarchy,
can be obtained as a reduction of
the Kadomtsev--Petviashvili (KP) hierarchy.
The latter is defined as the set of deformation equations
$${\partial L\over\partial t_k}=[(L^k)_+,L],$$
for the first order pseudo-differential operator
$$L\equiv
L(t,\partial)=\partial+u_1(t)\partial^{-1}+u_2(t)\partial^{-2}+\cdots,$$
here $\partial={\partial\over\partial t_1}$, $t=(t_1,t_2,\ldots)$ and
$(L^k)_+$ stands for the differential part of $L^k$.
Now $L$ dresses as $L=P\partial P^{-1}$ with
$$
P\equiv
P(t,\partial)=1+a_1(t)\partial^{-1}+a_2(t)\partial^{-2}+\cdots .$$
One can choose $P$ in such a way that
$$P(t,z)={{\tau (t-[z^{-1}])}\over {\tau (t)}},
$$
where $\tau (t)=\tau (t_1,t_2,t_3,\ldots) $ is the famous 
$\tau$--function, introduced by the Kyoto
group [DJKM1-3] and $[z]=(z,{z^2\over 2},{z^3\over 3},\ldots)$.
Sato [S] showed that such a $\tau$--function corresponds to a point of
some infinite Grassmannian $Gr$ (see e.g. [S,SW]). Let $H$ be the space of
formal Laurent series $\sum a_nt^n$, such that $a_n=0$ for $n>>0$. The
points of $Gr$ are those linear subspaces $W\subset H$ for which the
projection $\pi_+$ of $W$ into $H_+=\{\sum a_nt^n\in H | a_n=0
\text{ for all } n<0\}$ is a Fredholm operator. The $k$--th reduction 
or $k$--th Gelfand--Dickey hierarchy is obtained by assuming that 
$$L^k= (L^k)_+,$$
which corresponds to a $\tau$--function for which
$${{\partial \tau}\over{\partial t_k}}=\lambda\tau\qquad\text{for some
}\lambda\in \Bbb{C}.$$
In the polynomial case, i.e. $\tau$ is a polynomial, clearly
$\lambda=0$.
The point in the Grassmannian that corresponds to such a reduced 
$\tau$--function satisfies
$$t^kW\subset W.$$
In recent years a lot of attention has been drawn to a new kind of
reduction of the KP hierarchy, viz. the so--called 
$k$--constrained KP hierarchies [AFGZ,C,CWZ,CZ,D,DS,OS] 
(and references therein). Here one assumes that
$$L^k= (L^k)_+ +q\partial^{-1}r,\tag{1.1}$$
$q=q(t),r=r(t)$ being functions.
Under this condition the KP hierarchy is constrained to
$$
{\partial L\over\partial t_k}=[(L^k)_+,L],\qquad
{\partial q\over\partial t_k}=(L^k)_+, q\qquad
{\partial r\over\partial t_k}=-(L^k)^*_+ r.\tag{1.2}
$$
Here $A^*$ 
stands for the adjoined operator of $A$ (see e.g. [KV] for
more details about pseudo--differential operators). 
The AKNS, Yajima--Oikawa and Melnikov
hierarchies are some of the examples that appear amongst these
constrained KP families.

In this paper we consider the generalization of this $k$--constrained
KP hierarchy, wich was introduced by Sidorenko and Strampp in [SS], the
$n$--vector $k$--constrained hierarchy. We assume that
$$L^k= (L^k)_+ +\sum_{j=1}^n q_j\partial^{-1}r_j,\tag{1.3}$$
then one obtains the following integrable system:
$$
{\partial L\over\partial t_k}=[(L^k)_+,L],\qquad
{\partial q_j\over\partial t_k}=(L^k)_+ q_j,\qquad
{\partial r_j\over\partial t_k}=-(L^k)^*_+ r_j\qquad\text{for }1\le
j\le n.\tag{1.4}
$$
For $k=1$ this hierarchy contains the coupled vector non--linear
Schr\"odinger.
Zhang and Cheng showed in [ZC] that if one asumes that
$$q_j(t)={{\rho_j(t)}\over {\tau(t)}}\qquad\text{and}\qquad 
r_j(t)={{\sigma_j(t)}\over {\tau(t)}},\tag{1.5}$$
then $L$, $q_j$ and $r_j$, $1\le j\le n$ satisfy the $n$--vector
$k$--constrained hierarchy if and only if $\tau(t)$, $\rho_j(t)$ and
$\sigma_j(t)$ satisfy the following set of equations:
$$\align 
Res_{z=0}e^{-\eta(t,z)}\tau(t)e^{\xi(t,z)}
e^{\eta(t',z)}\tau(t')e^{-\xi(t',z)}
&=0,\tag{1.6}\\
Res_{z=0}z^ke^{-\eta(t,z)}\tau(t)e^{\xi(t,z)}
e^{\eta(t',z)}\tau(t')e^{-\xi(t',z)}
&=\sum_{j=1}^n \rho_j(t)\sigma_j(t'),\tag{1.7}\\
Res_{z=0}z^{-1}e^{-\eta(t,z)}\tau(t)e^{\xi(t,z)}
e^{\eta(t',z)}\rho_j(t')e^{-\xi(t',z)}
&=\rho_j(t)\tau(t'),\tag{1.8}\\
Res_{z=0}z^{-1}e^{-\eta(t,z)}\sigma_j(t)e^{\xi(t,z)}
e^{\eta(t',z)}\tau(t')e^{-\xi(t',z)}&=\tau(t)\sigma_j(t').\tag{1.9}
\endalign
$$
where
$$\eta(t,z)=\sum_{i=1}^\infty {{\partial}\over{\partial t_i}}
{{z^{-i}}\over i},\qquad
\xi(t,z)=\sum_{i=1}^\infty t_iz^i\tag{1.10}$$
and $Res_{z=0}\sum_i a_iz^i=a_{-1}$.

In the case that $n=1$, Loris and Willox [LW] show that one can deduce
some additional bilinear identities, but now involving ${{\partial
\tau}\over{\partial t_k}}$. It is unclear if this is possible for
$n>1$, but we will not need these extra bilinear identities.

We will show in this paper that in fact $L$ satisfies the $n$-vector
$k$--constrained KP hierarchy (1.3-4) if and only if
the corresponding point $W$ in $Gr$ has a linear
subspace $W'\subset W$ of codimension $n$ such that 
$$t^kW'\subset W.\tag{1.10}$$
We will proof this only in the polynomial case, i.e. 
polynomial $\tau$,$\rho_j$ and $\sigma_j$, but
we expect that this 
is still true in the non-polynomial case. We
use the representation theory of the infinite--dimensional matrix 
group $GL_\infty$, developed by Kac and Peterson [KP1-2] (see also
[KR]), to achieve this result.

Notice that in this way we get a filtration of hierarchies, i.e., the
$n$--vector $k$--constrained hierarchy is a subsystem of the 
$(n+1)$--vector $k$--constrained hierarchy, $n=0$ being the
$k$--reduced KP or Gelfand--Dickey hierarchies.
\vskip 10pt

\subheading{\S 2.  The semi-infinite wedge
representation of the group $GL_{\infty}$ and Sato's Grassmannian}

\vskip 10pt
\noindent Consider the infinite complex matrix group
$$GL_{\infty} = \{ A = (a_{ij})_{i,j \in {\Bbb Z}+\frac{1}{2}}|A\
\text{is invertible and all but a finite number of}\ a_{ij} -
\delta_{ij}\ \text{are}\ 0\}$$
and its Lie algebra
$$gl_{\infty} = \{ a = (a_{ij})_{i,j \in {\Bbb Z}+\frac{1}{2}}|\
\text{all but a finite number of}\ a_{ij}\ \text{are}\ 0\}$$
with bracket $[a,b] = ab-ba$.  The Lie algebra $gl_{\infty}$ has a
basis consisting of matrices $E_{ij},\ i,j \in {\Bbb Z} + \frac{1}{2}$, where
$E_{ij}$ is the matrix with a $1$ on the $(i,j)$-th entry and zeros
elsewhere.
Let ${\Bbb C}^{\infty} = \bigoplus_{j \in {\Bbb Z}+\frac{1}{2}} {\Bbb
C} v_{j}$ be an infinite dimensional complex vector space with fixed
basis $\{ v_{j}\}_{j \in {\Bbb Z}+\frac{1}{2}}$.  Both the group
$GL_{\infty}$ and its Lie algebra $gl_{\infty}$ act linearly on
${\Bbb C}^{\infty}$ via the usual formula:
$$E_{ij} (v_{k}) = \delta_{jk} v_{i}.$$

The well-known semi--infinite wedge representation is
constructed as follows [KP2] (see also [KR] and [KV]).  
The semi-infinite wedge space $F =
\Lambda^{\frac{1}{2}\infty} {\Bbb C}^{\infty}$ is the vector space
with a basis consisting of all semi-infinite monomials of the form
$v_{i_{1}} \wedge v_{i_{2}} \wedge v_{i_{3}} \ldots$, where $i_{1} >
i_{2} > i_{3} > \ldots$ and $i_{\ell +1} = i_{\ell} -1$ for $\ell >>
0$.  We can now define representations $R$ of $GL_{\infty}$ and $r$
of $gl_{\infty}$ on $F$ by
$$\aligned
R(A) (v_{i_{1}} \wedge v_{i_{2}} \wedge v_{i_{3}} \wedge \cdots) &= A
v_{i_{1}} \wedge Av_{i_{2}} \wedge Av_{i_{3}} \wedge \cdots ,\\ 
r(a) (v_{i_{1}} \wedge v_{i_{2}} \wedge v_{i_{3}} \wedge \cdots ) &=
\sum_{k} v_{i_{1}} \wedge v_{i_{2}} \wedge \cdots \wedge v_{i_{k-1}}
\wedge av_{i_{k}} \wedge v_{i_{k+1}} \wedge \cdots . 
\endaligned\tag{2.1}
$$
These equations are related by the usual formula:
$$\exp (r(a)) = R(\exp a)\ \text{for}\ a \in gl_{\infty}.$$
In order to perform calculations later on, it is convenient to
introduce a larger group
$$\align
\overline{GL_{\infty}} = \{ A = (a_{ij})_{i,j \in {\Bbb Z}+\frac{1}{2}}|&A\
\text{is invertible and all but a finite}\\
&\text{ number of}\ a_{ij} -
\delta_{ij}\ \text{with}\ i\ge j\ \text{are}\ 0\}
\endalign$$
and its Lie algebra
$$\overline{gl_{\infty}} = \{ a = (a_{ij})_{i,j \in {\Bbb Z}+\frac{1}{2}}|\
\text{all but a finite number of}\ a_{ij}\ \text{with}\ i\ge j\ \text{are}\ 0\}.$$
Both $\overline{GL_{\infty}}$ and $\overline{gl_{\infty}}$ act on a
completion $\overline{\Bbb C^\infty}$ of the space $\Bbb C^\infty$,
where 
$$\overline{\Bbb C^\infty}=\{\sum_j c_jv_j | c_j=0\ \text{for}\
j>>0\}.$$
It is easy to see that the representations $R$ and $r$ extend to
representations of $\overline{GL_{\infty}}$ and
$\overline{gl_{\infty}}
$ on the space $F$.

The representation $r$ of $gl_{\infty}$ and $\overline{gl_{\infty}}$ 
can be described
in  terms of wedging and contracting
operators in $F$
(see e.g. [KP2,KR]). Let $v_j^*$ be the linear functional on $\Bbb C^\infty$
defined by  $\langle v_i^*,v_j\rangle :=v_i^*(v_j)=\delta_{ij}$ and let $\Bbb C^{\infty *}=\bigoplus_{j \in {\Bbb Z}+\frac{1}{2}} {\Bbb
C} v_{j}^*$ be the restricted dual of $\Bbb C^{\infty}$,
then for any $w\in \Bbb C^{\infty }$, we define  a wedging operator $\psi^+(w)$
on $F$ by
$$\psi^{+}(w) (v_{i_{1}} \wedge v_{i_{2}} \wedge \cdots ) =
w\wedge v_{i_{1}} \wedge v_{i_{2}} \cdots .\tag{2.2}$$
Let $w^*\in\Bbb C^{\infty *}$, we define a contracting operator
$$\psi^{-}(w^*) (v_{i_{1}} \wedge v_{i_{2}} \wedge \cdots ) = 
\sum_{s=1}^\infty (-1)^{s+1} 
\langle w^*,v_{i_s}\rangle v_{i_{1}} \wedge v_{i_{2}} \wedge \cdots \wedge
v_{i_{s-1}} \wedge v_{i_{s+1}} \wedge \cdots .\tag {2.3}$$
For simplicity we write 
$$\psi^{+}_{j}=\psi^{+}(v_{-j}),\qquad
\psi^{-}_{j}=\psi^{-}(v_j^*)\qquad\text{for }j \in {\Bbb Z} +
\frac{1}{2}\tag{2.4}$$
These operators satisfy the following relations
$(i,j \in {\Bbb Z}+\frac{1}{2}, \lambda ,\mu = +,-)$:
$$\psi^{\lambda}_{i} \psi^{\mu}_{j} + \psi^{\mu}_{j}
\psi^{\lambda}_{i} = \delta_{\lambda ,-\mu} \delta_{i,-j},$$
hence they generate a Clifford algebra, which we denote by ${\Cal C}\ell$.

Introduce the following elements of $F$ $(m \in {\Bbb Z})$:
$$|m\rangle = v_{m-\frac{1}{2} } \wedge v_{m-\frac{3}{2} } \wedge
v_{m-\frac{5}{2} } \wedge \cdots .$$
It is clear that $F$ is an irreducible ${\Cal C}\ell$-module generated
by the vacuum $|0\rangle$ such that
$$\psi^{\pm}_{j} |0\rangle = 0 \ \text{for}\ j > 0 .$$
It is straightforward that the representation $r$ is given by the
following formula:
$$r(E_{ij}) = \psi^{+}_{-i} \psi^{-}_{j}. \tag{2.5}$$
Define the {\it charge decomposition}
$$F = \bigoplus_{m \in {\Bbb Z}} F^{(m)} \tag{2.6}$$
by letting
$$\text{charge}(|0\rangle ) = 0\ \text{and charge} (\psi^{\pm}_{j}) =
\pm 1. \tag{2.7}$$
It is clear that the charge decomposition is invariant with respect
to $r(g\ell_{\infty})$ (and hence with respect to $R(GL_{\infty})$).
Moreover, it is easy to see that each $F^{(m)}$ is irreducible with
respect to $g\ell_{\infty}$ (and $GL_{\infty}$).  Note that
$|m\rangle$ is its highest weight vector, i.e.
$$\align
&r(E_{ij})|m\rangle = 0 \ \text{for}\ i < j, \\
&r(E_{ii})|m\rangle = 0\  (\text{resp.}\ = |m\rangle ) \ \text{if}\ i > m\
(\text{resp. if}\ i < m).
\endalign
$$
\hskip 15pt
Let $w\in F$, we define the Annihilator space  $Ann (w)$ of $w$
as follows:
$$Ann(w)=\{v\in\Bbb C^\infty | v\wedge w=0\}.\tag{2.8}$$
Notice that $Ann(w)\ne 0$, since $v_j\in Ann(w)$ for $j<<0$.
This Annihilator space for perfect (semi--infinite) wedges 
$w\in F^{(m)}$ is related to the
$GL_{\infty}\text{-orbit}$
$${\Cal O}_m
= R(GL_{\infty})|m\rangle \subset F^{(m)}$$
of the highest weight vector $|m\rangle$ as follows.
Let $A=(A_{ij})_{i,j \in {\Bbb Z}}\in GL_\infty$, denote by 
$A_j=\sum_{i\in{\Bbb Z}} A_{ij}v_i$  then by (2.8)
$$\tau_m=R(A)|m\rangle =A_{m-{1\over 2}}\wedge
A_{m-{3\over 2}}\wedge A_{m-{5\over 2}}\wedge
\cdots,\tag{2.9}$$
with $A_{-j}=v_{-j}$ for $j>>0$.
Notice that since $\tau_m$ is a perfect
(semi-infinite) wedge
$$Ann(\tau_m)=\sum_{j<m} \Bbb CA_j
\subset \Bbb C^\infty.$$
By identifying $v_i=t^{-i-{1\over 2}}$ 
for $i\in\Bbb Z+{1\over 2}$, we can write 
$A_j=A_j(t)=\sum_{i\in\Bbb Z+{1\over 2}}
A_{ij}t^{-i-{1\over 2}}$ as a Laurent polynomial
in $t$. In this way we can identify
$Ann(\tau_m)$ with a subspace 
$W_{\tau_m}=\sum_{j<m}\Bbb C A_j(t)$ of the
space $H$ of all Laurent polynomials. Notice that this 
space $H$ differs from the one described in section 1.
So from now on let $Gr$ consist of all linear subspaces
of $H$ which contain 
$$H_j:=\sum_{i=-j}^\infty \Bbb C t^i$$
for $j>>0$ and let $Gr=\cup_{m\in\Bbb Z}Gr_m$ (disjoint union) with
$$Gr_m=\{W\in Gr| H_j\subset W\ \text{and}\ \dim W/H_j=m-j
\ \text{for}\ j<<0\},$$
then we can construct a cannonical map
$$\phi :{\Cal O}_m\to Gr_m,\qquad \phi(\tau_m)
=W_{\tau_m}:=\sum_{i<m} \Bbb C A_i(t).$$
It is clear that $\phi (|m\rangle )=H_{m}$ and that $\phi$ is 
surjective with fibers $\Bbb C^\times$. This construction 
is due to Sato [S], we call $Gr$ the polynomial Grassmannian.
From now on we will call a perfect wedge also a $\tau$--function (N.B.
$\tau=0$ is also a $\tau$--function).

\vskip 10pt
\subheading{\S 3. The boson-fermion correspondence}

\vskip 10pt
\noindent  
Introduce the  fermionic fields $(z \in {\Bbb C}^{\times})$:
$$\psi^{\pm }[z] \overset{\text{def}}\to{=} \sum_{k \in {\Bbb
Z}+\frac{1}{2}} \psi^{\pm
}_{k} z^{-k-\frac{1}{2}}.\tag{3.1}$$
Next we introduce bosonic fields:
$$\alpha[z] \equiv \sum_{k \in {\Bbb Z}} \alpha_{k} z^{-k-1}
\overset{def}\to{=} :\psi^{+}[z] \psi^{-}[z]:, \tag{3.2}$$
where $:\ :$ stands for the {\it normal ordered product} defined in
the usual way $(\lambda ,\mu = +$ or $-$):
$$:\psi^{\lambda }_{k} \psi^{\mu }_{\ell}: = \cases \psi^{
\lambda }_{k}
\psi^{\mu }_{\ell}\ &\text{if}\ \ell \ge k\\
-\psi^{\mu }_{\ell} \psi^{\lambda }_{k} &\text{if}\ \ell <k.
\endcases$$
One checks (using e.g.  Wick's formula) that the operators
$\alpha_{k}$ satisfy the commutation relations 
of the associative
oscillator algebra,  one has:
$$[\alpha_{k},\alpha_{\ell}] =
k\delta_{k,-\ell}\quad\text{and}
\quad\alpha_{k}|m\rangle = 0 \ \text{if}\ k > 0.\tag{3.3}$$
In order to express the fermionic fields $\psi^{\pm }(z)$ in terms of
the bosonic operators $\alpha_\ell$, we need some additional operator
$Q$.  This operators is uniquely defined as
follows:
$$Q(v_{i_{1}} \wedge v_{i_{2}} \wedge \cdots )
=(v_{i_{1}+1} \wedge v_{i_{2}+1} \wedge \cdots ).\tag{3.4}$$
So
$$Q|0\rangle = |1\rangle ,\ Q\psi^{\pm
}_{k} = \psi^{\pm}_{k\mp 1}Q$$
and $Q$ satisfies the following commutation relations with the $\alpha$'s:
$$[\alpha_{k},Q] =
\delta_{k0}Q.$$
In this paper  the operator $Q^{-k}$ will play an important role. If
$w_{m-{1\over 2}}\wedge w_{m-{3\over 2}}\wedge\cdots$ is a perfect
wedge then
$$Q^{-k}(w_{m-{1\over 2}}\wedge w_{m-{3\over 2}}\wedge\cdots)=
\Lambda^k w_{m-{1\over 2}}\wedge \Lambda^kw_{m-{3\over
2}}\wedge\cdots,\tag{3.5}$$
where $\Lambda=\sum_{j\in \Bbb Z+{1\over 2}}E_{j,j+1}$.
\proclaim{Theorem 3.1}{\bf ([DJKM1], [JM])}
$$\psi^{\pm}[z] = Q^{\pm 1}z^{\pm \alpha_{0}} \exp
(\mp \sum_{k < 0} \frac{1}{k} \alpha_{k}z^{-k})\exp(\mp
\sum_{k > 0} \frac{1}{k} \alpha_{k} z^{-k}). \tag{3.6}$$
\endproclaim
\demo{Proof} See [TV].
\enddemo
The operators on the right-hand side of (3.6) are called vertex
operators.  They made their first appearance in string theory (cf.
[FK]).
\vskip 10pt
\noindent We now describe the boson-fermion
correspondence.  Let ${\Bbb C}[t]$ be the space of polynomials in
indeterminates $t = (t_1,t_2,t_3,\ldots) $.  
Let $B ={\Bbb C}[q,q^{-1},t]= {\Bbb C}[t] \otimes_{\Bbb C} {\Bbb
C}[q,q^{-1}]$ 
be the tensor
product of algebras.  Then the boson-fermion
correspondence is the vector space isomorphism
$$\sigma :F @>\sim >> B,$$
given by
$$\sigma (\alpha_{-m_{1}} \ldots
\alpha_{-m_{s}}|k\rangle ) = m_{1} \ldots
m_{s}t_{m_{1}} \ldots t_{m_{s}} q^k.$$
Notice that the power of $q$ is the value of the charge.
The transported action of the operators $\alpha_{m}$ and $Q$ looks
as follows:
$$\sigma Q\sigma^{-1}=q,\qquad\sigma \alpha_m\sigma^{-1}=
\cases
-mt_m \ &\text{if}\ m <0, \\
{\partial\over {\partial t_m}}\ &\text{if}\ m >0, \\
q{\partial\over {\partial q}}\ &\text{if}\ m =0. 
\endcases \tag{3.7}$$
Hence 
$$\sigma \psi^\pm [z]\sigma^{-1}=q^{\pm 1}z^{\pm q{\partial\over {\partial
q}}}e^{\pm\xi(t,z)}
e^{\mp\eta(t,z)},\tag{3.8}$$
with
$\eta(t,z)$ and 
$\xi(t,z)$  given by (1.10)

\vskip 10pt
\subheading{\S 4. Identification of the bilinear identities}
\vskip 10pt
\noindent 
From now on we assume that $\tau\in F^{(m)}$, hence that $\tau$ is the
inverse immage under $\sigma$. Using the boson--fermion correspondence
of the previous section, we rewrite the bilinear identities (1.6-9) of
Zhang and Cheng now as equations in $F\otimes F$. Notice first the
following equality of operators on $F\otimes F$:
$$Res_{z=0} \psi^+[z]\otimes\psi^-[z]=
\sum_{i \in {\Bbb Z}+\frac{1}{2}} \psi^{+}_{i}\otimes \psi^{-}_{-i}$$
Now (1.6-9) turn into the following equations:
$$\align
\sum_{i \in {\Bbb Z}+\frac{1}{2}} \psi^{+}_{i} \tau \otimes \psi^{-}_{-i}
\tau&= 0, \tag{4.1}\\
\sum_{i \in {\Bbb Z}+\frac{1}{2}} \psi^{+}_{i} \tau \otimes \psi^{-}_{-i}
Q^{-k}\tau&=\sum_{j=1}^n \rho_j\otimes \sigma_j , \tag{4.2}\\
\sum_{i \in {\Bbb Z}+\frac{1}{2}} \psi^{+}_{i} \tau \otimes
\psi^{-}_{-i}\rho_j&=\rho_j\otimes\tau,\tag{4.3}\\
\sum_{i \in {\Bbb Z}+\frac{1}{2}} \psi^{+}_{i}\sigma_j\otimes \psi^{-}_{-i}
Q^{-k}\tau&=Q^{-k}\tau\otimes\sigma_j .\tag{4.4}
\endalign
$$
Here $Q^{-k}\tau\in F^{(m-k)}$, $\rho_j\in F^{(m+1)}$ and $\sigma_j\in
F^{(m-k-1)}$ for all $1\le j\le n$.
Equation (4.1) is called the KP hierarchy in the fermionic picture, it
characterizes
the $GL_\infty$-orbit ${\Cal O}_m$, i.e.
\proclaim{Proposition 4.1}{\bf ([KP2])} A non-zero element $\tau$ of $F^{(m)}$
lies in ${\Cal O_m}$ if and only if $\tau$ satisfies the equation $(4.1)$.
\endproclaim
If $\tau\in {\Cal O_m}$, then we can write $\tau$ as a perfect wedge
$$\tau=w_{m-{1\over 2}}\wedge w_{m-{3\over 2}}\wedge w_{m-
{5\over 2}}
\wedge w_{m-{7\over 2}}\wedge\cdots,\tag{4.5}$$
such that $w_{-\ell}=v_{-\ell}$ for $\ell>>0$. The corresponding point
$W_\tau\in Gr_m$ is then given by
$$W_\tau=\langle w_{m-{1\over 2}}, w_{m-{3\over 2}}, w_{m-
{5\over 2}},w_{m-{7\over 2}},\ldots\rangle.\tag{4.6}$$
The  geometrical interpretation of (4.3-4) is given by the following
proposition.
\proclaim {Proposition 4.2} Let $\tau\in{\Cal O}_m$, $\rho\in F^{(m+1)}$ and $\sigma\in
F^{(m-1)}$, then

\noindent {\bf (1)} $\tau$ and $\rho$ satisfy
$$\sum_{i \in {\Bbb Z}+\frac{1}{2}} \psi^{+}_{i} \tau \otimes
\psi^{-}_{-i}\rho=\rho\otimes\tau,\tag{4.7}$$
if and only if $\rho\in{\Cal O}_{m+1}$ and $W_\tau\subset W_\rho$,

\noindent {\bf (2)} $\tau$ and $\sigma$
satisfy
$$\sum_{i \in {\Bbb Z}+\frac{1}{2}} \psi^{+}_{i}\sigma\otimes \psi^{-}_{-i}
\tau=\tau\otimes\sigma ,\tag{4.8}$$ 
if and only if $\sigma\in{\Cal O}_{m-1}$ and
$W_\sigma\subset W_\tau$.
\endproclaim
\demo{Proof}
Without loss of generality we may assume (since the operator 
$\sum_i\psi^{+}_{i}\otimes \psi^{-}_{-i}$ commutes with the action of 
$R(GL_\infty)\otimes R(GL_\infty)$) that $\tau=|m\rangle$. Then (4.7)
is equivalent to
$$\sum_{i>m}v_i\wedge |m\rangle \otimes\psi^-_i\rho=\rho\otimes
|m\rangle.$$
Since all elements $v_i\wedge |m\rangle $, for $i>m$, are linearly
independent, we deduce that $\psi^-_i\rho=\lambda_i |m\rangle$ and
that
$\rho\in \langle v_i\wedge  |m\rangle|i>m\rangle$. 
Hence $\rho=w\wedge |m\rangle$ for some $w\in\Bbb C^\infty$ and thus
$\rho\in {\Cal O}_{m+1} $ and $W_\tau\subset W_\rho$.

The converse, since 
$W_\tau\subset W_\rho$, $\rho=w\wedge |m\rangle$ for some $w\in\Bbb
C^\infty$.
Then
$$
\aligned
\sum_{i \in {\Bbb Z}+\frac{1}{2}} \psi^{+}_{i} \tau \otimes
\psi^{-}_{-i}(w\wedge \tau)&=(w\wedge \tau)\otimes \tau 
-(1\otimes \psi^+(w))(
\sum_{i \in {\Bbb Z}+\frac{1}{2}} \psi^{+}_{i} \tau \otimes
\psi^{-}_{-i}\tau)\\
&=(w\wedge \tau)\otimes \tau
\endaligned 
$$
For $\tau=|m\rangle$, (4.8) is equivalent to
$$\sum_{i<m}(v_i\wedge\sigma)\otimes\psi_i^-|m\rangle=|m\rangle\otimes\sigma.
$$
Since the elements $\psi_i^-|m\rangle$ for $i<m$ are all linearly
independent, we conclude that $v_i\wedge \sigma=\lambda_i|m>$ and that
$\sigma\in \langle \psi^-_i|m\rangle |i<m\rangle$. Hence $\sigma=
\sum_{i=-\infty}^{m-{1\over 2}}a_i\psi^-_i|m\rangle$. Since $\sigma\in
F^{(m-1)}$, $a_i=0$ for all $i<-N<<0$. We now calculate $Ann(\sigma)$.
Clearly $Ann(\sigma)\subset \langle v_i|i<m\rangle= Ann (|m\rangle)$,
so let $v=\sum_{i<m}(-)^ib_iv_i$, then $\sum_{i=-N+{1\over
2}}^{m-{1\over 2}}a_ib_i=0$. Hence, if $\sigma\not =0$, we only find
one restriction for the collection of $b_i$'s, from which we conclude
that $\sigma$ is a perfect wedge. The converse of this statement
follows immediately by writing $\tau=w\wedge\sigma$.
\hfill $\square$
\enddemo
We next proof the following
\proclaim {Proposition 4.3} Let $\tau\in{\Cal O}_m$, 
$\rho_j\in{\Cal O}_{(m+1)}$ and $\sigma_j\in
{\Cal O}_{(m-k-1)}$, $1\le j\le n$, be related by
$$W_\tau\subset  W_{\rho_j},\qquad W_{\sigma_j}\subset
\Lambda^kW_\tau,\tag{4.9}$$
then $\tau$ satisfies equation $(4.2)$ if and only if there exists a
subspace $W'\subset W_\tau$ of codimension $n$ such that
$\Lambda^kW'\subset W_\tau$.
\endproclaim
\demo {Proof}
Notice first that $\Lambda^kW_\tau=W_{Q^{-k}\tau}$.
We assume that $n$ is minimal, so that all $\sigma_j$
and $\rho_j$ are nonzero perfect wedges, and that $\tau$ is of the form
(4.5). Then 
$$\align
\sum_{i \in {\Bbb Z}+\frac{1}{2}}& \psi^{+}_{i} \tau \otimes \psi^{-}_{-i}
Q^{-k}\tau\\
&=\sum_{\ell=0}^\infty (-)^\ell \Lambda^k w_{m-\ell-{1\over
2}}\wedge \tau\otimes \Lambda^kw_{m_{1\over 2}}\wedge\cdots \wedge\Lambda^k
w_{m-\ell+{1\over 2}}\wedge 
\Lambda^k w_{m-\ell-{3\over 2}}\wedge \cdots\\
&=\sum_{j=1}^n u_j\wedge \tau\otimes\sigma_j,
\endalign
$$
where
$\rho_j=u_j\wedge \tau$. Since all vectors
$\Lambda^kw_{m_{1\over 2}}\wedge\cdots\wedge \Lambda^k
w_{m-\ell+{1\over 2}}\wedge 
\Lambda^k w_{m-\ell-{3\over 2}}\wedge \cdots$ are linearly
independent, we deduce that 
$$\Lambda^k w_{m-\ell-{1\over
2}}\wedge u_1\wedge u_2\wedge\cdots\wedge u_n\wedge \tau=0,$$
for all $\ell=0,1,2,\dots$. Since we have assumed that $n$ is minimal,
also all $u_j$'s are linearly independent and moreover
$u_1\wedge u_2\wedge\cdots\wedge u_n\wedge \tau\not =0$,
hence
$$\Lambda^k w_{m-\ell-{1\over 2}}\in \langle
u_1,u_2,\ldots,u_n,w_{m-{1\over 2}}, w_{m-{3\over 2}},\ldots\rangle,$$
so there exists a subspace $W'\subset W_\tau$ of codimension $n$ such that
$\Lambda^kW'\subset W_\tau$.

For the converse, choose a basis $w_{m-n-{1\over 2}}, w_{m-n-{3\over
2}},\dots$ of $W'$ and extend it to a basis  
$w_{m-{1\over 2}}, w_{m-{3\over 2}},\ldots,w_{m-n+{1\over 2}}
w_{m-n-{1\over 2}}, w_{m-n-{3\over
2}},\dots$ of $W_\tau$, then
$$\align
\sum_{i \in {\Bbb Z}+\frac{1}{2}}& \psi^{+}_{i} \tau \otimes \psi^{-}_{-i}
Q^{-k}\tau\\
&=\sum_{\ell=0}^\infty (-)^\ell \Lambda^k w_{m-\ell-{1\over
2}}\wedge \tau\otimes \Lambda^kw_{m-{1\over 2}}\wedge\cdots\wedge \Lambda^k
w_{m-\ell+{1\over 2}}\wedge 
\Lambda^k w_{m-\ell-{3\over 2}}\wedge \cdots\\
&=\sum_{\ell=0}^{n-1} (-)^\ell \Lambda^k w_{m-\ell-{1\over
2}}\wedge \tau\otimes \Lambda^kw_{m-{1\over 2}}\wedge\cdots\wedge \Lambda^k
w_{m-\ell+{1\over 2}}\wedge 
\Lambda^k w_{m-\ell-{3\over 2}}\wedge \cdots.
\endalign
$$
So choose
$$\align
\rho_j&=\Lambda^k
w_{m-j+{1\over 2}}\wedge\tau,\\
\sigma_j&=\Lambda^kw_{m-{1\over 2}}\wedge\cdots \wedge\Lambda^k
w_{m-j+{3\over 2}}\wedge 
\Lambda^k w_{m-j-{1\over 2}}\wedge \cdots,
\endalign
$$
then $W_\tau,\ \Lambda^kW_{\tau},\ 
W_{\sigma_j}$ and $W_{\rho_j}$ clearly satisfy the equations
(4.9).\hfill $\square$
\enddemo
From this proposition we deduce the main Theorem of this paper
\proclaim {Theorem 4.4} The pseudo-differential operator
$$L=\partial+u_1\partial^{-1}+u_2\partial^{-2}+\cdots,$$
satisfies the $n$--vector $k$--constrained KP hierarchy if and only if
the corresponding point $W\in Gr_m$ has a subspace  $W'$ of
codimension $n$ such that $t^kW'\subset W$.
\endproclaim
As an easy consquence we obtain
\proclaim {Corollary 4.5} Let $\tau$ be a polynomial $\tau$--function of the 
$n$--vector $k$--constrained KP hierarchy, then 
${{\partial\tau}\over{\partial t_k}}=\sum_{\ell=1}^n \tau_\ell$
where every $\tau_\ell$ satisfies the KP hierarchy, i.e. equation
$(4.1)$.
\endproclaim
\demo {Proof} The proof follows immediately by taking the same basis
for $W_\tau$ as in the converse part of the proof of Proposition
4.3.\hfill $\square$
\enddemo
\noindent If $n=1$, one can proof [V] that every polynomial
$\tau$--function $\tau$, for which ${{\partial\tau}\over{\partial
t_k}}$ is again  $\tau$--function, is a solution  of the
$k$--constrained KP hierarchy.

Notice that we have constructed a natural filtration  on the space
$Gr_m$, which is determined by the $n$--vector $k$--constrained KP
hierarchy for $n=0,1,2,\ldots$. Let 
$$\aligned
Gr_m^{(n,k)}=\{W\in Gr_m|&\text{there exists a subspace } W'\subset W\\
&\text{of codimension } n\ \text{such that } t^kW'\subset
W\},\endaligned
\tag{4.10}$$
then 
$$Gr_m^{(0,k)}\subset Gr_m^{(1,k)}\subset\cdots\subset 
Gr_m^{(n,k)}\subset Gr_m^{(n+1,k)}\subset\cdots .\tag{4.11}$$
It is obvious that every point $W\in Gr_m$ (in this polynomial case)
is contained in $Gr_m^{(n,k)}$ for  $n>>0$, in other words
$$Gr_m=\bigcup_{n\in\Bbb Z_+}Gr_m^{(n,k)}.\tag{4.12}$$
So for every $\tau$--function of the KP hierarchy there exists a
non--negative integer $n$ such that for all $m\ge n$, $\tau$ is also
a $\tau$--function of the $m$--vector $k$--constrained KP hierarchy.
In other words, for every $L$, corresponding to a polynomial
$\tau$--function, one can find a non--negative integer $n$ such that
$L$ satisfies (1.3).
\vskip 10pt
\subheading{\S 5. Polynomial solutions of the $n$--vector
$k$--constrained KP hierarchy.}
\vskip 10pt
\noindent 
We will now state an immediate consequence of the boson--fermion
correspondence, viz., 
we calculate the immage under $\sigma$
of a perfect wedge of the form (2.9). One finds the
folowing result
\proclaim{Proposition 5.1}Let $S_i$ be the elementary Schur functions,
defined by $\exp \sum_{i=1}^\infty t_iz^i=\sum_{i\in\Bbb Z}^\infty
S_i(t)z^i$ ($S_i=0$ for $i<0$) and let $\tau_m\in{\Cal O}_m$
be of the form $(2.9)$, i.e.,
$$\tau_m=A_{m-{1\over 2}}\wedge
A_{m-{3\over 2}}\wedge A_{m-{5\over 2}}\wedge\cdots$$
with $A_j=\sum_{i\in\Bbb Z+{1\over2}} A_{ij} v_i$ and
$A_{-k}=v_{-k}$ for all $k>N>>0$. Set 
$A= (A_{ij})_{i\in\Bbb Z+{1\over 2},m>j\in\Bbb Z+{1\over 2}}$
and let $\Lambda=\sum_{i\in \Bbb Z={1\over 2}} E_{i,i+1}
\in \overline{gl_\infty}$. Then
$$\sigma(\tau_m)=det \left (\sum_{i,j=-n+{1\over 2}}^{m-{1\over
2}}\left (\sum_{\ell=-N+{1\over 2}}^\infty S_{\ell-i}A_{\ell
j}\right )E_{ij}\right )q^m.\tag{5.1}$$
\endproclaim
\demo {Proof}
The proof of this proposition is the same as the proof of
Theorem 6.1 of [KR]. One computes 
$$\sigma (\exp\left (\sum_{i=1}^\infty t_i\Lambda^i\right )\tau_m)$$
and takes the coefficient of $q^m$. One thus obtains (see also [DJKM1,M]):
$$\sigma(\tau_m)=
\det \left (\left (\exp\left (\sum_{i=1}^\infty t_i\Lambda^i\right
)A\right )_{<m}\right )q^m,\tag{5.2}$$
where $B_{<m}$ denotes
the submatrix of $B$ where one only takes the rows 
$j\in\Bbb Z+{1\over 2}$ with $j<m$.
Notice that $\sum_i t_i\Lambda^i
\in\overline{gl_\infty}$ and $\exp(\sum_i t_i\Lambda^i)
\in\overline{GL_\infty}$. Here we calculate the determinant 
of an infinite matrix, however there is no problem,
since the matrix is of the form $(B_{ij})_{m>i,j\in\Bbb Z
+{1\over 2}}$ with all but a finite number of 
$B_{ij}-\delta_{ij}$ with $i\ge j$ are zero.

It is clear that one can subtract $\sum_{i<-N}A_{ij}v_i$ from every $A_j$, with
$j>-N$, in $\tau_m$,
this will not change $\tau_m$. Then the new $A$ is of the form
$$A=\sum_{-N<i,-N<j<m}A_{ij} E_{ij}+\sum_{i<-N}E_{ii},$$ 
it is then
straightforward, using the elementary Schur functions,  
to calculate the right-hand-side of (5.2). One finds
formula (5.1).\hfill$\square$
\enddemo
We will use this proposition to obtain all polynomial solutions of the
$n$--vector $k$--constrained KP hierarchy. Notice that our approach is
different from the one in [ZC]. Instead of taking $\tau_m$
of the form (2.9), we may choose another basis of $W_{\tau_m}$ and
construct the corresponding perfect wedge, it is clear that this wil
be a multiple of $\tau_m$. We can choose this basis in such a way
$$W_{\tau_m}=\langle A_{m-{1\over 2}},
A_{m-{3\over 2}},A_{m-{5\over 2}},\ldots, A_{-N+{1\over 2}}, 
v_{-N-{1\over 2}},v_{-N-{3\over 2}},\ldots\rangle,$$
such that $A_j=\sum_{i=-N+{1\over 2}}^\infty A_{ij}v_i$ and that, 
except for at most $n$ vectors $A_j$, all  $A_j$ satisfy
the following condition
$$
\Lambda^kA_j\cases
=A_\ell \ \text{for some }-N+{1\over 2}\le \ell\le m-{1\over 2},\text{ or}\\
\in\langle v_{-N-{1\over 2}},v_{-N-{3\over
2}},\ldots\rangle.\endcases
$$
Of course every $A_j$ is bounded, i.e., there exists an integer $M$ such
that all $A_j=\sum_{i=-N+{1\over 2}}^{M-{1\over 2}}A_{ij}v_i$.
Now making a shift in the index and permuting the columns 
we obtain the following result:
\proclaim {Proposition 5.2} Let $M,N\in \Bbb Z$ such that $M>N>0$ and let $e_j$, $1\le j\le M$ be an orthonormal
basis of $\Bbb C^M$. Let $R$ be the $M\times M$--matrix
$R=\sum_{i=1}^{M-k}E_{i,i+k}$ and let $A=(A_{ij})_{1\le i\le M,1\le j\le N}$ be an
$M\times N$--matrix of rank $N$. Denote by
$A_j=\sum_{i=1}^MA_{ij}e_i$.
If all $A_j$ satisfy the condition that 
$RA_j\ne A_i$ for all $1\le i<j$ 
and if all $A_j$, except for at most $n$,
satisfy the condition that
$$
RA_j=\cases
A_{j+1}\qquad\text{or}\\
0,\endcases
$$
then 
$$\tau=det \left (\sum_{i,j=1}^N\left (\sum_{\ell=1}^M S_{\ell-i}A_{\ell
j}\right )E_{ij}\right )\tag{5.3}$$
is a $\tau$--function of the $n$--vector $k$--constrained KP
hierarchy.
All polynomial solutions can be obtained in this way.
\endproclaim 
\vskip 10pt
\noindent {\bf Acknowledgements} It is a pleasure to thank
Gerard Helminck, Ignace Loris and Gerhard Post for helpful discussions, 
and Walter Strampp for drawing my attention to this subject.

\enddocument